\documentclass[11pt,a4paper]{amsart}
\usepackage[T1]{fontenc}
\usepackage[utf8]{inputenc}
\usepackage{lmodern}
\usepackage[a4paper,textwidth=15.6cm,textheight=23.4cm,centering,headsep=0.8cm]{geometry}
\usepackage{amsmath,amssymb,amsthm,mathtools}
\usepackage{booktabs,array}
\usepackage[final]{microtype}
\usepackage{enumitem}
\usepackage{xcolor}
\usepackage{tikz}
\usetikzlibrary{arrows.meta,calc,fit,backgrounds}
\usepackage{hyperref}
\hypersetup{
  pdfauthor={Santiago Pineda Montoya, Edisson Gallego Gonz\'alez, Jaiberth Porras Barrera, Johan Heriberto R\'ua Mu\~noz},
  pdftitle={Laplacian eigenvectors of generically weighted graphs have pairwise distinct coordinates},
  pdfsubject={Spectral graph theory},
  pdfkeywords={weighted graph Laplacian; Fiedler vector; generic weights; structural genericity; Krylov determinant},
  colorlinks=true,linkcolor=blue!45!black,citecolor=blue!45!black,urlcolor=blue!45!black}
\allowdisplaybreaks[2]
\numberwithin{equation}{section}
\setlist[enumerate,1]{label=\textup{(\roman*)},leftmargin=*,widest=iii,itemsep=2pt,topsep=4pt}

\newtheorem{theorem}{Theorem}[section]
\newtheorem{lemma}[theorem]{Lemma}
\newtheorem{proposition}[theorem]{Proposition}
\newtheorem{corollary}[theorem]{Corollary}
\theoremstyle{definition}
\newtheorem{definition}[theorem]{Definition}
\newtheorem{example}[theorem]{Example}
\theoremstyle{remark}
\newtheorem{remark}[theorem]{Remark}
\newtheorem{question}[theorem]{Question}

\newcommand{\one}{\mathbf 1}
\newcommand{\R}{\mathbb R}
\newcommand{\Z}{\mathbb Z}
\newcommand{\Q}{\mathbb Q}
\newcommand{\cS}{\mathcal S}
\newcommand{\cZ}{\mathcal Z}
\newcommand{\cP}{\mathcal P}
\newcommand{\spec}{\operatorname{spec}}

\begin{document}

\title[Distinct coordinates of generic Laplacian eigenvectors]{Laplacian eigenvectors of generically weighted graphs\\have pairwise distinct coordinates}
\author[S. Pineda Montoya]{Santiago Pineda Montoya}
\address{Universidad Nacional de Colombia}
\email{sapinedamo@unal.edu.co}
\thanks{Corresponding author: Santiago Pineda Montoya (\texttt{sapinedamo@unal.edu.co}).}

\author[E. Gallego Gonz\'alez]{Edisson Gallego Gonz\'alez}
\address{Universidad EAFIT}
\email{egalleg8@eafit.edu.co}

\author[J. Porras Barrera]{Jaiberth Porras Barrera}
\address{Universidad Nacional de Colombia}
\email{jporras@unal.edu.co}

\author[J. H. R\'ua Mu\~noz]{Johan Heriberto R\'ua Mu\~noz}
\address{Universidad de Antioquia}
\email{heriberto.rua@udea.edu.co}
\date{16 November 2026}
\subjclass[2020]{Primary 05C50; Secondary 15A18, 05C22, 47A55}
\keywords{Weighted graph Laplacian, Fiedler vector, algebraic connectivity, generic weights, structural genericity, Krylov determinant}

\begin{abstract}
Poignard, Pereira and Pade conjectured that, for a generic choice of the nonzero weights of a given graph, the Fiedler vector of the weighted Laplacian has pairwise distinct coordinates. We give a proof that also covers every nonconstant eigenvector. For each connected graph $\Gamma$ on $n$ vertices we construct a nonzero homogeneous polynomial $P_\Gamma$ with integer coefficients in the edge weights, of degree $n(n-1)(n^2-n+2)/4$, which vanishes at a positive weighting exactly when the Laplacian has a repeated eigenvalue, an eigenvector with a zero coordinate, or a nonconstant eigenvector with two equal coordinates. The exceptional weights therefore lie in an algebraic hypersurface, which gives structural genericity in the sense of Poignard, Pereira and Pade; algebraically independent weights are never exceptional, and there are good integer weights not exceeding the degree plus one. The polynomial is a product of Krylov determinants. Its nonvanishing follows by induction on trees: when a leaf is attached with a small weight, a cut identity gives the first-order profile of the new eigenvector and separates the coordinates that coincide in the limit. Exact small cases show that exceptional weightings need not be symmetric.
\end{abstract}

\maketitle

\section{Introduction}

\subsection{The conjecture}
Let $\Gamma=(V,E)$ be a finite simple graph with vertex set $V=\{1,\ldots,n\}$, and let $w=(w_e)_{e\in E}$ be positive edge weights. The weighted Laplacian is
\begin{equation}\label{eq:laplacian}
 L_\Gamma(w)=\sum_{\{i,j\}\in E}w_{ij}\,(e_i-e_j)(e_i-e_j)^T ,
\end{equation}
that is, $L_\Gamma(w)=D-W$ with $W$ the weighted adjacency matrix and $D$ the diagonal matrix of weighted degrees. It is symmetric and positive semidefinite, and it annihilates the constant vector $\one$. When $\Gamma$ is connected its eigenvalues can be written
\[
 0=\lambda_1<\lambda_2\le\cdots\le\lambda_n .
\]
The number $\lambda_2$ is the algebraic connectivity of the weighted graph~\cite{Fiedler73}, and an eigenvector of $\lambda_2$ is a \emph{Fiedler vector}. Its coordinates order and split the vertices in spectral partitioning and seriation~\cite{Fiedler75,ABH}, and they govern how diffusively coupled networks gain or lose synchronization when the coupling structure is modified~\cite{PPP}.

Poignard, Pereira and Pade~\cite{PPP} studied which spectral properties hold for a generic choice of the weights when the graph, that is, the set of nonzero weights, is kept fixed. For connected undirected graphs with positive weights they proved that the following properties are generic in this structural sense: all eigenvalues are simple~\cite[Theorem~3.1]{PPP}; the Fiedler vector has no zero entry~\cite[Theorem~3.2]{PPP}; and there is a basis of eigenvectors without zero entries~\cite[Corollary~5.1]{PPP}. Their paper closes with a conjecture~\cite[Section~7]{PPP}: for a generic choice of the nonzero weights of a given graph, ``the Fiedler vector has only distinct coordinates''. We show that the same holds for every nonconstant eigenvector simultaneously and that the exceptional weights are cut out by a single explicit polynomial; simplicity of the spectrum and the absence of zero coordinates come out of the same argument.

\subsection{Results}
Throughout, an \emph{eigenvector} is a nonzero vector, and all statements about eigenvectors are invariant under rescaling.

\begin{definition}[Separating weights]\label{def:separating}
Let $\Gamma$ be a connected graph on $n\ge2$ vertices. A weight vector $w\in(0,\infty)^E$ is \emph{separating} if
\begin{enumerate}[label=\textup{(S\arabic*)},leftmargin=3.2em]
\item\label{S1} the $n$ eigenvalues of $L_\Gamma(w)$ are pairwise distinct;
\item\label{S2} no eigenvector of $L_\Gamma(w)$ has a zero coordinate;
\item\label{S3} every eigenvector of $L_\Gamma(w)$ orthogonal to $\one$ has pairwise distinct coordinates.
\end{enumerate}
We write $\cS_\Gamma\subset(0,\infty)^E$ for the set of separating weights.
\end{definition}

Under~\ref{S1} the eigenvectors orthogonal to $\one$ are exactly the eigenvectors of $\lambda_2,\ldots,\lambda_n$, that is, the nonconstant eigenvectors. In particular, for separating weights the Fiedler vector is unique up to a scalar, and its coordinates are nonzero and pairwise distinct.

\begin{theorem}[Separating weights are generic]\label{thm:main}
Let $\Gamma=(V,E)$ be a connected graph on $n\ge2$ vertices. There is a nonzero homogeneous polynomial $P_\Gamma\in\Z[w_e:e\in E]$ of degree
\begin{equation}\label{eq:degree}
 d_n=\frac{n(n-1)(n^2-n+2)}4
\end{equation}
such that, for every $w\in(0,\infty)^E$,
\[
 w\ \text{is separating}\quad\Longleftrightarrow\quad P_\Gamma(w)\ne0 .
\]
Consequently $\cS_\Gamma$ is an open dense cone in $(0,\infty)^E$, and its complement in $(0,\infty)^E$ is a relatively closed set of Lebesgue measure zero, contained in the real algebraic set $\{P_\Gamma=0\}$.
\end{theorem}

The polynomial is explicit: $P_\Gamma(w)=\cP_n(L_\Gamma(w))$, where $\cP_n$ is the product of Krylov determinants defined in~\eqref{eq:Pn}. The theorem is the equivalence of Proposition~\ref{prop:certificate}, which is linear algebra, together with the assertion that $P_\Gamma$ is not the zero polynomial, which requires one separating weighting of one spanning tree and is the substance of the proof.

In the language of~\cite{PPP}, Theorem~\ref{thm:main} reads as follows.

\begin{corollary}[Structural genericity; the conjecture]\label{cor:ppp}
Let $\Gamma$ be a connected graph on $n\ge2$ vertices.
\begin{enumerate}
\item\label{ppp:dense} For every $w\in(0,\infty)^E$ and every $\epsilon>0$ there is $w'\in\cS_\Gamma$ with $\max_{e\in E}|w'_e-w_e|<\epsilon$.
\item\label{ppp:null} For every $w\in(0,\infty)^E$, the set of perturbations $t\in\R^E$ for which $w+t$ lies in $(0,\infty)^E$ and is not separating has Lebesgue measure zero.
\end{enumerate}
Hence the property ``$\lambda_2$ is simple and the Fiedler vector has pairwise distinct coordinates'', and more strongly the conjunction of \ref{S1}--\ref{S3}, is structurally generic among the Laplacian matrices of connected undirected weighted graphs in the sense of~\cite[Definition~2.2]{PPP}. In particular the conjecture of~\cite[Section~7]{PPP} holds.
\end{corollary}

Since the polynomial has integer coefficients, genericity also takes arithmetic and probabilistic forms.

\begin{corollary}[Arithmetic and probabilistic forms]\label{cor:arithmetic}
Let $\Gamma$ be a connected graph on $n\ge2$ vertices and $d_n$ the degree~\eqref{eq:degree}.
\begin{enumerate}
\item If the positive numbers $(w_e)_{e\in E}$ are algebraically independent over $\Q$, then $w$ is separating.
\item If $w$ is a random vector whose law is absolutely continuous with respect to Lebesgue measure on $(0,\infty)^E$, then $w$ is separating almost surely.
\item If the weights are independent and uniformly distributed on $\{1,\ldots,N\}$, then $w$ fails to be separating with probability at most $d_n/N$. In particular there are separating integer weights with $1\le w_e\le d_n+1$ for every edge.
\end{enumerate}
\end{corollary}

The examples of Section~\ref{sec:examples} show what the exceptional set looks like. Symmetric weightings are always exceptional when $n\ge3$ (Proposition~\ref{prop:symmetry}), which explains why unweighted graphs so often have Fiedler vectors with repeated entries; the converse fails, since some weightings without any symmetry have a Fiedler vector with a zero coordinate or with two equal coordinates (Examples~\ref{ex:path4} and~\ref{ex:cycle4}). For the path on three vertices and for the triangle the exceptional sets are computed exactly (Examples~\ref{ex:path3} and~\ref{ex:triangle}), and on every weighted path the Fiedler vector is strictly monotone, although it may vanish at a vertex and the other eigenvectors may have ties (Proposition~\ref{prop:paths}).

\subsection{Strategy}
The proof combines an algebraic certificate with an explicit witness.

\emph{A certificate.} For a symmetric matrix $A$ with $A\one=0$, the Krylov matrix of the coordinate vector $e_v$ has nonzero determinant exactly when $A$ has simple spectrum and no eigenvector vanishes at $v$; the nonvanishing of this determinant is Kalman's rank condition for controllability~\cite{Kalman}. In the same way, the Krylov matrix of the difference $e_i-e_j$ inside the invariant hyperplane $\one^\perp$ has full rank exactly when the spectrum on $\one^\perp$ is simple and no eigenvector there takes equal values at $i$ and $j$. Since the columns of this matrix lie in $\one^\perp$, one maximal minor already detects its rank. Multiplying the determinants for all vertices and these minors for all pairs gives a polynomial $\cP_n(A)$ in the entries of $A$ whose nonvanishing is equivalent to \ref{S1}--\ref{S3} (Proposition~\ref{prop:certificate}).

\emph{A witness.} The certificate reduces genericity to a single separating weighting, which shows that $P_\Gamma=\cP_n\circ L_\Gamma$ is not identically zero. Setting to zero the weights outside a spanning tree reduces this to trees (Lemma~\ref{lem:integrality}), and trees are handled by induction, attaching one leaf $u$ with a small weight $\varepsilon$ to a tree $T$ that already has separating weights. The eigenvectors that continue those of $T$ keep their separated coordinates and acquire a small nonzero coordinate at $u$. The difficulty is the new eigenvector, whose eigenvalue $\mu(\varepsilon)$ tends to zero: in the limit it equals $1$ on all of $T$ and $-n$ at $u$, so its coordinates on $T$ coincide to order zero. A summation identity over the two sides of an edge (Lemma~\ref{lem:cut}) gives the expansion
\[
 x_v-x_a=\mu\,h_v+O(\mu^2),\qquad h_v=\sum_{e\in[a,v]}\frac{|D_e|}{w_e},
\]
where $a$ is the vertex of attachment, $[a,v]$ is the path from $a$ to $v$ and $D_e$ is the set of vertices cut off from $a$ by the edge $e$. The error is quadratic in $\mu$, with an explicit constant (Lemma~\ref{lem:neweigenvector}(iii)), and this precision is what allows the first-order term to separate the coordinates. The numbers $h_v$ are pairwise distinct for all weights outside a nowhere dense set, and then the coordinates are separated for small $\varepsilon$ (Proposition~\ref{prop:leaf}).

\subsection{Related work}
Genericity of simple spectrum, and of nondegeneracy properties of eigenfunctions, is classical for Laplace operators on manifolds under perturbations of the metric~\cite{Uhlenbeck}. For graphs with a fixed support the corresponding statements are due to Poignard, Pereira and Pade~\cite{PPP}. The argument given here is of a different kind: a polynomial identity takes the place of the perturbation of an arbitrary weighting, so that a single good weighting implies genericity. In this way it also yields the three genericity results of~\cite{PPP} quoted above, for undirected graphs, together with the separation of coordinates.

Krylov determinants of graph matrices are the subject of controllability theory. Godsil asked how often the pair formed by an adjacency matrix and the all-ones vector is controllable~\cite{Godsil}, and O'Rourke and Touri proved that this holds for almost all graphs~\cite{ORT}; simplicity of the spectrum of random symmetric matrices was established by Tao and Vu~\cite{TaoVu}. Those are statements about random discrete matrices. Here the support is arbitrary and fixed, and the weights are continuous parameters, which makes an algebraic argument available. Nodal properties of Laplacian eigenvectors on graphs are surveyed in~\cite{BLS}.

\subsection{Organization}
Section~\ref{sec:prelim} collects notation, a lemma on zero sets of polynomials and two elementary perturbation facts. Section~\ref{sec:certificate} constructs the polynomial $\cP_n$ and proves the certificate. Section~\ref{sec:trees} treats trees. Section~\ref{sec:proofs} proves Theorem~\ref{thm:main} and its corollaries. Section~\ref{sec:examples} contains exact examples, the relation with symmetry and the case of paths, and Section~\ref{sec:questions} collects remarks and open questions.

\section{Preliminaries}\label{sec:prelim}

\subsection{Notation}
Vectors in $\R^n$ are columns indexed by $V=\{1,\ldots,n\}$; $e_1,\ldots,e_n$ is the standard basis, $\one=e_1+\cdots+e_n$, and $x_v$ or $x(v)$ denotes the coordinate of $x$ at $v$. The Euclidean norm and the induced operator norm are both written $\|\cdot\|$. We write
\[
 \cZ_n=\{A\in\R^{n\times n}:A=A^T,\ A\one=0\}
\]
for the space of symmetric matrices with zero row sums. For $A\in\cZ_n$ the hyperplane $\one^\perp$ is invariant, since $\langle Ax,\one\rangle=\langle x,A\one\rangle=0$. Conditions \ref{S1}--\ref{S3} make sense for every $A\in\cZ_n$ in place of $L_\Gamma(w)$, and we call $A$ \emph{separating} when they hold.

A \emph{graph} is finite, simple and undirected. For a graph $\Gamma=(V,E)$ and an arbitrary real vector $w\in\R^E$, formula~\eqref{eq:laplacian} defines $L_\Gamma(w)\in\cZ_n$; the map $w\mapsto L_\Gamma(w)$ is linear. If $E'\subset E$ and $w$ vanishes outside $E'$, then
\begin{equation}\label{eq:restriction}
 L_\Gamma(w)=L_{\Gamma'}(w|_{E'}),\qquad\Gamma'=(V,E').
\end{equation}

\begin{lemma}\label{lem:psd}
Let $\Gamma$ be a graph and $w\in(0,\infty)^E$. Then
\[
 x^TL_\Gamma(w)x=\sum_{\{i,j\}\in E}w_{ij}(x_i-x_j)^2\qquad(x\in\R^n),
\]
so $L_\Gamma(w)$ is positive semidefinite, and its kernel consists of the vectors that are constant on each connected component of $\Gamma$. If $\Gamma$ is connected, the kernel is spanned by $\one$, every nonzero eigenvalue is positive, and its eigenvectors are orthogonal to $\one$.
\end{lemma}
\begin{proof}
The identity is immediate from~\eqref{eq:laplacian}. A positive semidefinite matrix annihilates exactly the vectors on which its quadratic form vanishes, which here are the vectors with $x_i=x_j$ for every edge. Eigenvectors of a symmetric matrix for different eigenvalues are orthogonal.
\end{proof}

\subsection{Zero sets of polynomials}

\begin{lemma}\label{lem:zeroset}
Let $P\in\R[x_1,\ldots,x_m]$ be a nonzero polynomial and $Z(P)=\{x\in\R^m:P(x)=0\}$. Then $Z(P)$ is closed, has Lebesgue measure zero, and has empty interior. Consequently, for every nonempty open set $U\subset\R^m$, the set $U\setminus Z(P)$ is open, dense in $U$, and of full measure in $U$.
\end{lemma}
\begin{proof}
$Z(P)$ is closed because $P$ is continuous. We prove by induction on $m$ that it is a null set. For $m=1$ it is finite. Let $m>1$ and write $P=\sum_{j=0}^{d}P_j(x')\,x_m^j$ with $x'=(x_1,\ldots,x_{m-1})$ and $P_d\ne0$. By induction the set $N'=\{x':P_d(x')=0\}$ is null in $\R^{m-1}$. For $x'\notin N'$ the polynomial $P(x',\cdot)$ is nonzero, so the slice $\{x_m:P(x',x_m)=0\}$ is finite. By Tonelli's theorem, applied to the indicator function of the closed set $Z(P)$, the measure of $Z(P)$ is the integral over $x'$ of the one-dimensional measures of its slices, which vanish outside the null set $N'$. A null set has empty interior, and the remaining assertions follow.
\end{proof}

\subsection{Two perturbation facts}
For a real symmetric matrix $A$ of order $N$ let $\lambda_1(A)\le\cdots\le\lambda_N(A)$ be its eigenvalues. Weyl's inequality~\cite[Section~4.3]{HornJohnson} states that
\begin{equation}\label{eq:weyl}
 |\lambda_k(A)-\lambda_k(B)|\le\|A-B\|\qquad(1\le k\le N)
\end{equation}
for real symmetric $A$ and $B$.

\begin{lemma}[Stability of a simple eigenpair]\label{lem:stability}
Let $A$ be real symmetric with a simple eigenvalue $\lambda$, unit eigenvector $\phi$, and gap
$\delta=\min\{|\lambda-\lambda'|:\lambda'\in\spec(A),\ \lambda'\ne\lambda\}>0$. Let $B$ be real symmetric with $\rho=\|B-A\|<\delta/2$. Then $B$ has exactly one eigenvalue $\lambda_B$, counted with multiplicity, in the interval $(\lambda-\delta/2,\lambda+\delta/2)$; it satisfies $|\lambda_B-\lambda|\le\rho$, and it has a unit eigenvector $\phi_B$ with
\[
 \|\phi_B-\phi\|\le\frac{2\sqrt2\,\rho}{\delta}.
\]
\end{lemma}
\begin{proof}
Let $\lambda=\lambda_k(A)$. By~\eqref{eq:weyl}, $\lambda_B:=\lambda_k(B)$ satisfies $|\lambda_B-\lambda|\le\rho<\delta/2$, while for $j\ne k$ we have $|\lambda_j(B)-\lambda|\ge|\lambda_j(A)-\lambda|-\rho\ge\delta-\rho>\delta/2$. This proves the first two assertions. Let $y$ be a unit eigenvector of $B$ for $\lambda_B$ and write $y=c\,\phi+r$ with $r\perp\phi$. Then
\[
 (A-\lambda)r=(A-\lambda)y=(A-B)y+(\lambda_B-\lambda)y ,
\]
whose norm is at most $2\rho$. On $\phi^\perp$ the symmetric operator $A-\lambda$ has all its eigenvalues of absolute value at least $\delta$, so $\|(A-\lambda)r\|\ge\delta\|r\|$ and $\|r\|\le2\rho/\delta<1$. Replacing $y$ by $-y$ if necessary we may assume $c\ge0$, so $c=\sqrt{1-\|r\|^2}$ and
\[
 \|y-\phi\|^2=(1-c)^2+\|r\|^2=2-2\sqrt{1-\|r\|^2}\le2\|r\|^2 ,
\]
because $1-\sqrt{1-s}\le s$ for $0\le s\le1$. Take $\phi_B=y$.
\end{proof}

\section{A polynomial certificate}\label{sec:certificate}

\subsection{Krylov matrices}
For a square matrix $A$, a vector $b$ and an integer $r\ge1$, the \emph{Krylov matrix} is the matrix with columns $b,Ab,\ldots,A^{r-1}b$.

\begin{lemma}[Krylov determinants]\label{lem:krylov}
Let $A$ be a real symmetric $n\times n$ matrix and $H\subset\R^n$ an $A$-invariant subspace of dimension $r\ge1$. Let $u_1,\ldots,u_r$ be an orthonormal basis of $H$ consisting of eigenvectors of $A$, with $Au_k=\mu_ku_k$. For $b\in H$ put $K=[\,b,Ab,\ldots,A^{r-1}b\,]\in\R^{n\times r}$. Then
\begin{equation}\label{eq:gramkrylov}
 \det(K^TK)=\prod_{k=1}^{r}\langle u_k,b\rangle^2\prod_{1\le k<l\le r}(\mu_l-\mu_k)^2 .
\end{equation}
If $H=\R^n$, so that $r=n$ and $K$ is square, then
\begin{equation}\label{eq:krylovdet}
 \det K=\pm\prod_{k=1}^{n}\langle u_k,b\rangle\prod_{1\le k<l\le n}(\mu_l-\mu_k).
\end{equation}
\end{lemma}
\begin{proof}
Since $b=\sum_k\langle u_k,b\rangle u_k$, we have $A^jb=\sum_k\langle u_k,b\rangle\mu_k^ju_k$ for $j\ge0$. Hence
\[
 K=U\,C\,M,\qquad U=[\,u_1,\ldots,u_r\,],\quad C=\operatorname{diag}\bigl(\langle u_k,b\rangle\bigr)_{k=1}^r,\quad M=\bigl(\mu_k^{\,j-1}\bigr)_{k,j=1}^{r}.
\]
The matrix $M$ is a Vandermonde matrix, with $\det M=\prod_{k<l}(\mu_l-\mu_k)$. Since $U^TU=I_r$, we get $K^TK=M^TC^2M$, which gives~\eqref{eq:gramkrylov}. If $r=n$, then $U$ is orthogonal, $\det U=\pm1$, and $\det K=\det U\det C\det M$.
\end{proof}

In words: $\det(K^TK)\ne0$ if and only if the restriction of $A$ to $H$ has $r$ distinct eigenvalues and $b$ is orthogonal to none of its eigenvectors. Such a vector $b$ is called \emph{cyclic} for $A|_H$.

\begin{lemma}[Reduced Gram determinants]\label{lem:reduced}
Let $K\in\R^{n\times(n-1)}$ be a matrix whose columns lie in $\one^\perp$, and let $\tilde K$ be the square matrix obtained by deleting its last row. Then
\[
 \det(K^TK)=n\,(\det\tilde K)^2 .
\]
\end{lemma}
\begin{proof}
For $r\in V$ let $K^{(r)}$ be the matrix obtained from $K$ by deleting row $r$, so that $\tilde K=K^{(n)}$. By the Cauchy--Binet formula,
\[
 \det(K^TK)=\sum_{r=1}^{n}\bigl(\det K^{(r)}\bigr)^2 .
\]
Let $\rho_1,\ldots,\rho_n$ be the rows of $K$. Since every column lies in $\one^\perp$, $\rho_1+\cdots+\rho_n=0$. Let $r<n$. Adding to the row $\rho_r$ of $K^{(n)}$ all its other rows does not change the determinant and replaces $\rho_r$ by $\rho_1+\cdots+\rho_{n-1}=-\rho_n$. The resulting matrix differs from $K^{(r)}$ only by the sign of one row and the position of that row. Hence $|\det K^{(r)}|=|\det\tilde K|$ for every $r$, and the sum equals $n(\det\tilde K)^2$.
\end{proof}

\subsection{The polynomial}
Let $n\ge2$. For $A\in\cZ_n$, a vertex $v$ and a pair $i<j$ define
\begin{align}
 \kappa_v(A)&=\det\bigl[\,e_v,Ae_v,\ldots,A^{n-1}e_v\,\bigr],\label{eq:kappa}\\
 K_{ij}(A)&=\bigl[\,e_i-e_j,A(e_i-e_j),\ldots,A^{n-2}(e_i-e_j)\,\bigr]\in\R^{n\times(n-1)},\notag\\
 \delta_{ij}(A)&=\det\tilde K_{ij}(A),\label{eq:delta}
\end{align}
where $\tilde K_{ij}(A)$ is obtained from $K_{ij}(A)$ by deleting its last row, and
\begin{equation}\label{eq:Pn}
 \cP_n(A)=\prod_{v=1}^{n}\kappa_v(A)\cdot\prod_{1\le i<j\le n}\delta_{ij}(A),\qquad
 P_\Gamma(w)=\cP_n\bigl(L_\Gamma(w)\bigr)\quad(w\in\R^E).
\end{equation}
The columns of $K_{ij}(A)$ lie in $\one^\perp$, because $e_i-e_j\in\one^\perp$ and $\one^\perp$ is $A$-invariant; by Lemma~\ref{lem:reduced}, $\det\bigl(K_{ij}(A)^TK_{ij}(A)\bigr)=n\,\delta_{ij}(A)^2$. For $n=2$ the matrix $\tilde K_{12}(A)$ is the $1\times1$ matrix $(1)$ and $\delta_{12}=1$.

\begin{proposition}[Certificate]\label{prop:certificate}
Let $n\ge2$ and $A\in\cZ_n$. Let $u_1=\one/\sqrt n,u_2,\ldots,u_n$ be an orthonormal basis of eigenvectors of $A$, with eigenvalues $\mu_1=0,\mu_2,\ldots,\mu_n$. Then
\begin{align}
 \kappa_v(A)&=\pm\frac1{\sqrt n}\prod_{k=2}^{n}u_k(v)\prod_{1\le k<l\le n}(\mu_l-\mu_k),\label{eq:kappaformula}\\
 \delta_{ij}(A)^2&=\frac1n\prod_{k=2}^{n}\bigl(u_k(i)-u_k(j)\bigr)^2\prod_{2\le k<l\le n}(\mu_l-\mu_k)^2 .\label{eq:deltaformula}
\end{align}
Consequently $\cP_n(A)\ne0$ if and only if $A$ is separating, that is, if and only if \ref{S1}, \ref{S2} and \ref{S3} hold for $A$.
\end{proposition}
\begin{proof}
Such a basis exists because $\one$ is an eigenvector for the eigenvalue $0$ and $\one^\perp$ is invariant. Formula~\eqref{eq:kappaformula} is~\eqref{eq:krylovdet} with $b=e_v$, since $\langle u_k,e_v\rangle=u_k(v)$ and $u_1(v)=n^{-1/2}$. Formula~\eqref{eq:deltaformula} follows from~\eqref{eq:gramkrylov}, applied with $H=\one^\perp$, $r=n-1$, the basis $u_2,\ldots,u_n$ and $b=e_i-e_j\in\one^\perp$, since $\langle u_k,e_i-e_j\rangle=u_k(i)-u_k(j)$, together with $\det(K_{ij}^TK_{ij})=n\,\delta_{ij}^2$.

Suppose $\cP_n(A)\ne0$. Then $\kappa_1(A)\ne0$, so by~\eqref{eq:kappaformula} the numbers $\mu_1,\ldots,\mu_n$ are pairwise distinct, which is~\ref{S1}. Every eigenvector is then a multiple of some $u_k$. The eigenvector $u_1$ has no zero coordinate, and $u_k(v)\ne0$ for $k\ge2$ and every $v$ because $\kappa_v(A)\ne0$; this is~\ref{S2}. The eigenvectors orthogonal to $\one$ are the multiples of $u_2,\ldots,u_n$, and $u_k(i)\ne u_k(j)$ for $i<j$ because $\delta_{ij}(A)\ne0$; this is~\ref{S3}. Conversely, if \ref{S1}--\ref{S3} hold, every factor on the right of~\eqref{eq:kappaformula} and~\eqref{eq:deltaformula} is nonzero.
\end{proof}

\begin{lemma}[Integrality, degree, restriction]\label{lem:integrality}
Let $\Gamma=(V,E)$ be a graph on $n\ge2$ vertices.
\begin{enumerate}
\item $P_\Gamma$ is a polynomial in the variables $(w_e)_{e\in E}$ with integer coefficients. It is either zero or homogeneous of degree $d_n=n(n-1)(n^2-n+2)/4$.
\item If $\Gamma'=(V,E')$ with $E'\subset E$, and $w\in\R^E$ vanishes outside $E'$, then $P_\Gamma(w)=P_{\Gamma'}(w|_{E'})$.
\item If $\Gamma$ is disconnected, then $P_\Gamma=0$.
\end{enumerate}
\end{lemma}
\begin{proof}
(i) The entries of $L=L_\Gamma(w)$ are linear forms in $w$ with integer coefficients, so the entries of $L^k$ are integer polynomials, each zero or homogeneous of degree $k$. In the expansion of the determinant~\eqref{eq:kappa} every term takes one entry from each of the columns $e_v,Le_v,\ldots,L^{n-1}e_v$; hence $\kappa_v$ is zero or homogeneous of degree $0+1+\cdots+(n-1)=n(n-1)/2$. In the same way, the column of $\tilde K_{ij}$ formed by the first $n-1$ entries of $L^k(e_i-e_j)$ has degree $k$ for $0\le k\le n-2$, so $\delta_{ij}$ is zero or homogeneous of degree $0+1+\cdots+(n-2)=(n-1)(n-2)/2$. Adding the degrees of the $n$ factors $\kappa_v$ and the $\binom n2$ factors $\delta_{ij}$ gives
\[
 n\cdot\frac{n(n-1)}2+\binom n2\frac{(n-1)(n-2)}2=\frac{n(n-1)}2\cdot\frac{2n+(n-1)(n-2)}2=d_n .
\]
(ii) follows from~\eqref{eq:restriction}, because $P_\Gamma(w)$ depends on $w$ only through the matrix $L_\Gamma(w)$. (iii) If $\Gamma$ is disconnected, then for every $w\in\R^E$ the matrix $L_\Gamma(w)$ annihilates the indicator vector of each connected component, so $0$ is a multiple eigenvalue and \ref{S1} fails; by Proposition~\ref{prop:certificate}, $P_\Gamma(w)=0$ for all $w$.
\end{proof}

\begin{remark}
The factors $\kappa_v$ alone certify \ref{S1} and \ref{S2}, which are the properties shown to be structurally generic in~\cite{PPP}. The factors $\delta_{ij}$ are responsible for~\ref{S3}: by Lemma~\ref{lem:reduced}, $\delta_{ij}(A)$ vanishes exactly when $e_i-e_j$ fails to be cyclic for the restriction of $A$ to $\one^\perp$. The Gram determinants $\det(K_{ij}^TK_{ij})=n\,\delta_{ij}^2$ have the same zero sets, but using them would double the contribution of these factors to the degree. The degree $d_n$ is not claimed to be minimal; see Question~\ref{q:degree}.
\end{remark}

\section{Trees}\label{sec:trees}
We now show that every tree on at least two vertices has separating weights (Proposition~\ref{prop:trees}).

\subsection{A summation identity}

\begin{lemma}[Cut identity]\label{lem:cut}
Let $G$ be a graph with real edge weights and Laplacian $L$, and let $D$ be a set of vertices such that exactly one edge $e=\{p,q\}$, with $q\in D$ and $p\notin D$, joins $D$ to its complement. If $Lx=\mu x$, then
\[
 w_e\,(x_q-x_p)=\mu\sum_{v\in D}x_v .
\]
\end{lemma}
\begin{proof}
By~\eqref{eq:laplacian}, $(Lx)_v=\sum_{v'}w_{vv'}(x_v-x_{v'})$, the sum running over the neighbours $v'$ of $v$. Sum this over $v\in D$. An edge with both ends in $D$ contributes $w_{vv'}(x_v-x_{v'})+w_{vv'}(x_{v'}-x_v)=0$, and the only other edge meeting $D$ is $e$, which contributes $w_e(x_q-x_p)$. The left side of $Lx=\mu x$ therefore sums to $w_e(x_q-x_p)$ and the right side to $\mu\sum_{v\in D}x_v$.
\end{proof}

\subsection{Attaching a leaf}
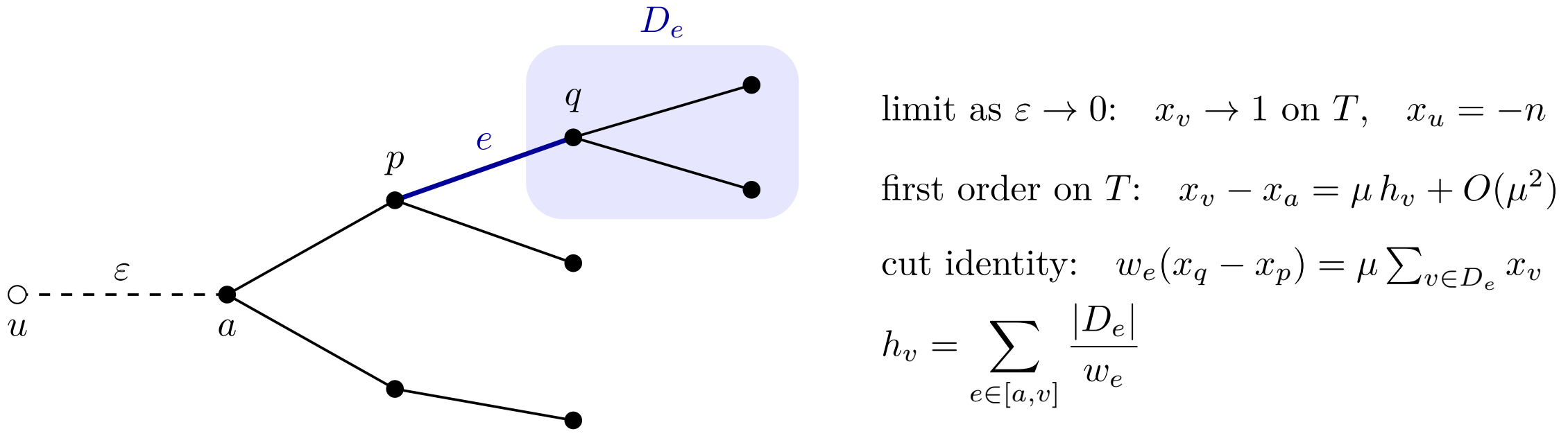
\begin{figure}[tbp]
\centering
\begin{tikzpicture}[x=1cm,y=1cm,font=\small,
  vtx/.style={circle,fill=black,inner sep=1.7pt},
  edge/.style={line width=0.7pt},
  new/.style={line width=0.7pt,dashed}]
  \coordinate (a) at (0,0);
  \coordinate (u) at (-2.0,0);
  \coordinate (p) at (1.6,0.9);
  \coordinate (r) at (1.6,-0.9);
  \coordinate (q) at (3.3,1.5);
  \coordinate (s) at (3.3,0.3);
  \coordinate (q1) at (5.0,2.0);
  \coordinate (q2) at (5.0,1.0);
  \coordinate (r1) at (3.3,-1.2);
  \begin{scope}[on background layer]
    \fill[blue!10,rounded corners=10pt] (2.85,0.72) rectangle (5.45,2.38);
  \end{scope}
  \draw[edge] (a)--(p) (a)--(r) (p)--(s) (q)--(q1) (q)--(q2) (r)--(r1);
  \draw[edge,blue!60!black,line width=1.3pt] (p)--(q) node[midway,above=1pt] {$e$};
  \draw[new] (a)--(u) node[midway,above] {$\varepsilon$};
  \foreach \x in {a,p,r,q,s,q1,q2,r1} \node[vtx] at (\x) {};
  \node[circle,draw=black,fill=white,inner sep=1.7pt] at (u) {};
  \node[below=3pt] at (a) {$a$};
  \node[above=3pt] at (p) {$p$};
  \node[above=3pt] at (q) {$q$};
  \node[below=3pt] at (u) {$u$};
  \node[blue!60!black] at (4.15,2.6) {$D_e$};
  \node[anchor=west] at (6.1,1.75) {limit as $\varepsilon\to0$:\quad $x_v\to1$ on $T$,\quad $x_u=-n$};
  \node[anchor=west] at (6.1,1.0) {first order on $T$:\quad $x_v-x_a=\mu\,h_v+O(\mu^2)$};
  \node[anchor=west] at (6.1,0.25) {cut identity:\quad $w_e(x_q-x_p)=\mu\sum_{v\in D_e}x_v$};
  \node[anchor=west] at (6.1,-0.6) {$h_v=\displaystyle\sum_{e\in[a,v]}\frac{|D_e|}{w_e}$};
\end{tikzpicture}
\caption{Attaching a leaf $u$ at the vertex $a$ of a tree $T$ with a small weight $\varepsilon$. The edge $e=\{p,q\}$ separates the shaded vertex set $D_e$ from $a$. The eigenvector of the new small eigenvalue $\mu=\mu(\varepsilon)$ is constant on $T$ in the limit; the cut identity gives its first-order profile $\mu h_v$, which separates the coordinates when the numbers $h_v$ are pairwise distinct.}
\label{fig:leaf}
\end{figure}

Fix the following data for this subsection (Figure~\ref{fig:leaf}).
\begin{itemize}[leftmargin=*,itemsep=2pt]
\item $T$ is a tree with vertex set $V=\{1,\ldots,n\}$, $n\ge2$, and $w\in\cS_T$ is a separating weight vector. Let $L=L_T(w)$, with eigenvalues $0=\lambda_1<\lambda_2<\cdots<\lambda_n$ and unit eigenvectors $\phi_1=\one/\sqrt n,\phi_2,\ldots,\phi_n$.
\item $a\in V$ is a vertex, and $T'$ is the tree obtained from $T$ by adding a new vertex $u$ and the edge $\{a,u\}$. For $\varepsilon>0$ we give this edge the weight $\varepsilon$ and write $L_\varepsilon=L_{T'}(w,\varepsilon)$. Ordering the vertices of $T'$ as $1,\ldots,n,u$,
\begin{equation}\label{eq:Leps}
 L_\varepsilon=\begin{pmatrix}L&0\\0&0\end{pmatrix}
 +\varepsilon\,(e_a-e_u)(e_a-e_u)^T .
\end{equation}
\item For an edge $e$ of $T$, the graph $T-e$ has two connected components; $D_e\subset V$ denotes the vertex set of the component that does not contain $a$. For $v\in V$, $[a,v]$ denotes the set of edges of the path from $a$ to $v$ in $T$, and
\begin{equation}\label{eq:h}
 h_v=\sum_{e\in[a,v]}\frac{|D_e|}{w_e},\qquad h_a=0 .
\end{equation}
\end{itemize}
We use the positive constants
\[
 \delta=\min_{1\le k<n}(\lambda_{k+1}-\lambda_k),\qquad
 \eta=\min\Bigl\{|\phi_k(v)|,\ |\phi_k(i)-\phi_k(j)|:\ 2\le k\le n,\ v\in V,\ i\ne j\Bigr\},
\]
which are positive because $w$ is separating, and
\[
 c=\frac{(n-1)^2}{\min_ew_e},\qquad h_{\max}=\max_{v\in V}h_v>0 .
\]

\begin{lemma}[Spectrum of $L_\varepsilon$]\label{lem:leafspectrum}
Let $0<\varepsilon<\delta/4$. Then $L_\varepsilon$ has $n+1$ distinct eigenvalues: the eigenvalue $0$; one eigenvalue $\mu(\varepsilon)\in(0,2\varepsilon]$; and, for each $k=2,\ldots,n$, exactly one eigenvalue $\lambda_k(\varepsilon)$ with $|\lambda_k(\varepsilon)-\lambda_k|\le2\varepsilon$.
\end{lemma}
\begin{proof}
By~\eqref{eq:Leps}, $L_\varepsilon-L_0=\varepsilon(e_a-e_u)(e_a-e_u)^T$ has norm $2\varepsilon$, and the eigenvalues of $L_0=L\oplus0$ in nondecreasing order are $0,0,\lambda_2,\ldots,\lambda_n$. Let $\nu_1\le\cdots\le\nu_{n+1}$ be the eigenvalues of $L_\varepsilon$. By~\eqref{eq:weyl}, $|\nu_1|,|\nu_2|\le2\varepsilon$ and $|\nu_{k+1}-\lambda_k|\le2\varepsilon$ for $2\le k\le n$. Since $T'$ is connected and its weights are positive, Lemma~\ref{lem:psd} gives $\nu_1=0<\nu_2$; put $\mu(\varepsilon)=\nu_2$ and $\lambda_k(\varepsilon)=\nu_{k+1}$. Because $\lambda_2\ge\delta$ and consecutive $\lambda_k$ differ by at least $\delta>4\varepsilon$, the intervals $[0,2\varepsilon]$ and $[\lambda_k-2\varepsilon,\lambda_k+2\varepsilon]$, $2\le k\le n$, are pairwise disjoint. Hence the $n+1$ eigenvalues are distinct.
\end{proof}

\begin{lemma}[Eigenvectors continuing those of $T$]\label{lem:oldeigenvectors}
Let $0<\varepsilon<\min\{\delta/4,\ \delta\eta/(12\sqrt2)\}$ and $2\le k\le n$. Every eigenvector of $L_\varepsilon$ for $\lambda_k(\varepsilon)$ has nonzero, pairwise distinct coordinates.
\end{lemma}
\begin{proof}
In the spectrum $\{0,\lambda_2,\ldots,\lambda_n\}$ of $L_0$ the eigenvalue $\lambda_k$ is simple, with unit eigenvector $(\phi_k,0)$ and gap at least $\delta$. Apply Lemma~\ref{lem:stability} with $A=L_0$, $B=L_\varepsilon$ and $\rho=2\varepsilon<\delta/2$: the eigenvalue of $L_\varepsilon$ in $(\lambda_k-\delta/2,\lambda_k+\delta/2)$ is $\lambda_k(\varepsilon)$, and it has a unit eigenvector $x$ with
\[
 \|x-(\phi_k,0)\|\le\frac{4\sqrt2\,\varepsilon}{\delta}<\frac\eta3 .
\]
Every coordinate of $x$ on $V$ is therefore within $\eta/3$ of the corresponding coordinate of $\phi_k$. By the definition of $\eta$,
\[
 |x_v|>\tfrac23\eta\quad(v\in V),\qquad |x_i-x_j|>\tfrac13\eta\quad(i\ne j\text{ in }V).
\]
The row of $u$ in $L_\varepsilon x=\lambda_k(\varepsilon)x$ reads $\varepsilon(x_u-x_a)=\lambda_k(\varepsilon)x_u$, that is,
\[
 x_u=-\frac{\varepsilon\,x_a}{\lambda_k(\varepsilon)-\varepsilon},\qquad
 \lambda_k(\varepsilon)-\varepsilon\ge\lambda_k-3\varepsilon\ge\delta-3\varepsilon>\frac\delta4 .
\]
Thus $x_u\ne0$, because $x_a\ne0$, and $|x_u|<4\varepsilon/\delta<\eta/3<|x_v|$ for every $v\in V$, so $x_u$ differs from all other coordinates. Since $\lambda_k(\varepsilon)$ is simple, every eigenvector is a multiple of $x$.
\end{proof}

\begin{lemma}[The new eigenvector]\label{lem:neweigenvector}
Let $0<\varepsilon<\min\{\delta/4,\ 1/(4c)\}$, let $\mu=\mu(\varepsilon)$, and let $x$ be an eigenvector of $L_\varepsilon$ for $\mu$.
\begin{enumerate}
\item $\sum_{v\in V}x_v\ne0$. Normalizing $x$ so that $\sum_{v\in V}x_v=n$, we have $x_u=-n$.
\item With this normalization, $|x_v-1|\le2c\mu$ for every $v\in V$; in particular $x_v>0$.
\item With this normalization, for all $v,v'\in V$,
\[
 \bigl|\,x_v-x_{v'}-\mu\,(h_v-h_{v'})\,\bigr|\le4c\,h_{\max}\,\mu^2 .
\]
\end{enumerate}
\end{lemma}
\begin{proof}
(i) Since $\mu\ne0$, Lemma~\ref{lem:psd} applied to $T'$ gives $x\perp\one$, that is, $\sum_{v\in V}x_v+x_u=0$. Suppose $\sum_{v\in V}x_v=0$. Then $x_u=0$, and the row of $u$, $\varepsilon(x_u-x_a)=\mu x_u$, gives $x_a=0$. By~\eqref{eq:Leps}, the rows of $V$ then read $Ly=\mu y$ for $y=(x_v)_{v\in V}$, and $y\ne0$ because $x\ne0$. So $\mu$ would be an eigenvalue of $L$ in $(0,2\varepsilon]\subset(0,\delta)$, whereas $L$ has no eigenvalue in $(0,\delta)$. This contradiction proves $\sum_{v\in V}x_v\ne0$, and the normalization gives $x_u=-n$.

(ii) Let $e=\{p,q\}$ be an edge of $T$ with $q\in D_e$. The vertex $a$ does not belong to $D_e$, so no vertex of $D_e$ is adjacent to $u$, and $e$ is the only edge of $T'$ joining $D_e$ to its complement. Lemma~\ref{lem:cut} gives
\begin{equation}\label{eq:cutT}
 x_q-x_p=\frac\mu{w_e}\sum_{v\in D_e}x_v .
\end{equation}
Let $M=\max_{v\in V}|x_v|$. Since $|D_e|\le n-1$, the right side of~\eqref{eq:cutT} is at most $\mu(n-1)M/\min_ew_e$ in absolute value. Two vertices of $V$ are joined in $T$ by a path with at most $n-1$ edges, so $|x_v-x_{v'}|\le c\mu M$ for all $v,v'\in V$. The mean of the $x_{v'}$ over $V$ equals $1$; hence $|x_v-1|\le c\mu M$ and $M\le1+c\mu M$. As $c\mu\le2c\varepsilon\le1/2$, we get $M\le2$ and $|x_v-1|\le2c\mu<1$.

(iii) Let $a=v_0,v_1,\ldots,v_r=v$ be the path from $a$ to $v$ and $e_s=\{v_{s-1},v_s\}$. The end of $e_s$ in $D_{e_s}$ is $v_s$, so summing~\eqref{eq:cutT} over $s$ gives
\[
 x_v-x_a=\mu\sum_{e\in[a,v]}\frac1{w_e}\sum_{v'\in D_e}x_{v'} .
\]
By (ii), $\bigl|\sum_{v'\in D_e}x_{v'}-|D_e|\bigr|\le2c\mu|D_e|$. Comparing with~\eqref{eq:h},
\[
 |x_v-x_a-\mu h_v|\le2c\mu^2h_v\le2c\,h_{\max}\,\mu^2 ,
\]
and the claim follows by subtracting the same estimate for $v'$.
\end{proof}

\begin{proposition}[Leaf attachment]\label{prop:leaf}
In the setting above, suppose that the numbers $h_v$, $v\in V$, are pairwise distinct, and let $\theta=\min_{v\ne v'}|h_v-h_{v'}|>0$. Then $(w,\varepsilon)\in\cS_{T'}$ for every
\[
 0<\varepsilon<\varepsilon_0:=\min\Bigl\{\frac\delta4,\ \frac{\delta\eta}{12\sqrt2},\ \frac1{4c},\ \frac{\theta}{8c\,h_{\max}}\Bigr\}.
\]
\end{proposition}
\begin{proof}
By Lemma~\ref{lem:leafspectrum} the spectrum of $L_\varepsilon$ is simple, which is~\ref{S1}, and every eigenvector is a multiple of $\one$, of an eigenvector of some $\lambda_k(\varepsilon)$, or of an eigenvector of $\mu(\varepsilon)$. The first has no zero coordinate. The second kind has nonzero, pairwise distinct coordinates by Lemma~\ref{lem:oldeigenvectors}. For the third, normalize as in Lemma~\ref{lem:neweigenvector}: then $x_v>0$ on $V$ and $x_u=-n<0$, so no coordinate vanishes and $x_u$ differs from every $x_v$. For $v\ne v'$ in $V$, part (iii) of that lemma and $\mu\le2\varepsilon$ give
\[
 |x_v-x_{v'}|\ge\mu\,\bigl(\theta-4c\,h_{\max}\,\mu\bigr)\ge\mu\,\bigl(\theta-8c\,h_{\max}\,\varepsilon\bigr)>0 .
\]
Hence \ref{S2} and \ref{S3} hold.
\end{proof}

\subsection{Every tree has separating weights}

\begin{lemma}[The profile $h$ is generically injective]\label{lem:hgeneric}
Let $T$ be a tree on $V$, $n\ge2$, let $a\in V$, and let $v\ne v'$ be vertices. The set of $w\in(0,\infty)^{E(T)}$ with $h_v(w)=h_{v'}(w)$ is closed in $(0,\infty)^{E(T)}$ and has empty interior.
\end{lemma}
\begin{proof}
The paths from $a$ to $v$ and from $a$ to $v'$ share an initial segment and then separate; the symmetric difference of $[a,v]$ and $[a,v']$ is the edge set of the path from $v$ to $v'$, which is nonempty. By~\eqref{eq:h},
\[
 h_v-h_{v'}=\sum_{e\in[a,v]\setminus[a,v']}\frac{|D_e|}{w_e}-\sum_{e\in[a,v']\setminus[a,v]}\frac{|D_e|}{w_e},
\]
with at least one term and all $|D_e|\ge1$. Multiplying by $\prod_ew_e$ gives a polynomial $Q_{vv'}$ in which the monomials $\prod_{e'\ne e}w_{e'}$, for the edges $e$ of the path from $v$ to $v'$, are distinct and have the nonzero coefficients $\pm|D_e|$. Thus $Q_{vv'}\ne0$, and on $(0,\infty)^{E(T)}$ the equation $h_v=h_{v'}$ is equivalent to $Q_{vv'}(w)=0$. Lemma~\ref{lem:zeroset} applies.
\end{proof}

\begin{proposition}[Trees]\label{prop:trees}
For every tree $T$ on $n\ge2$ vertices, $\cS_T\ne\varnothing$.
\end{proposition}
\begin{proof}
By induction on $n$. For $n=2$ the Laplacian is $w\left(\begin{smallmatrix}1&-1\\-1&1\end{smallmatrix}\right)$, with eigenvalues $0$ and $2w$ and eigenvectors $(1,1)^T$ and $(1,-1)^T$, so every $w>0$ is separating.

Let $T'$ be a tree on $n+1\ge3$ vertices. It has a leaf $u$; let $a$ be its neighbour and $T=T'-u$, a tree on $n$ vertices. By the induction hypothesis $\cS_T\ne\varnothing$, and by Proposition~\ref{prop:certificate}, $\cS_T=\{w\in(0,\infty)^{E(T)}:P_T(w)\ne0\}$ is open. By Lemma~\ref{lem:hgeneric}, the set of weights for which two of the numbers $h_v$ coincide is the zero set in $(0,\infty)^{E(T)}$ of the nonzero polynomial $\prod_{v<v'}Q_{vv'}$, with $Q_{vv'}$ as in the proof of that lemma, and so has empty interior by Lemma~\ref{lem:zeroset}. It therefore does not contain the nonempty open set $\cS_T$. Choose $w\in\cS_T$ with pairwise distinct $h_v$. By Proposition~\ref{prop:leaf}, $(w,\varepsilon)\in\cS_{T'}$ for all sufficiently small $\varepsilon>0$.
\end{proof}

\section{Proofs of the main results}\label{sec:proofs}

\begin{proof}[Proof of Theorem~\ref{thm:main}]
Let $w\in(0,\infty)^E$. By Proposition~\ref{prop:certificate} applied to $A=L_\Gamma(w)$, the weight vector $w$ is separating if and only if $P_\Gamma(w)=\cP_n(L_\Gamma(w))\ne0$.

We show that $P_\Gamma$ is not the zero polynomial. Choose a spanning tree $T$ of $\Gamma$ and, by Proposition~\ref{prop:trees}, a weight vector $w_T\in\cS_T$. Let $w^\circ\in\R^E$ agree with $w_T$ on the edges of $T$ and vanish on the other edges. By Lemma~\ref{lem:integrality}(ii) and Proposition~\ref{prop:certificate},
\[
 P_\Gamma(w^\circ)=P_T(w_T)\ne0 .
\]
By Lemma~\ref{lem:integrality}(i), $P_\Gamma$ has integer coefficients and is homogeneous of degree $d_n$.

Finally, $\cS_\Gamma=(0,\infty)^E\setminus Z(P_\Gamma)$. By Lemma~\ref{lem:zeroset} it is open and dense in $(0,\infty)^E$, and its complement there is relatively closed and null. It is a cone because $P_\Gamma$ is homogeneous.
\end{proof}

When $\Gamma$ is not a tree, the witness $w^\circ$ lies on the boundary of the positive orthant, and no perturbation is needed to move it inside: a polynomial that is nonzero at one point is nonzero, and a nonzero polynomial cannot vanish on the open set $(0,\infty)^E$.

\begin{proof}[Proof of Corollary~\ref{cor:ppp}]
\ref{ppp:dense} is the density of $\cS_\Gamma$ in $(0,\infty)^E$. For \ref{ppp:null}, the set in question is $\bigl((0,\infty)^E\setminus\cS_\Gamma\bigr)-w$, a translate of a null set.

It remains to match these statements with the definitions of~\cite{PPP}. There, a property of matrices in a class $\mathcal M\subset\cZ_n$ is called structurally generic when two conditions hold~\cite[Definitions~2.1 and~2.2]{PPP}. First, every $M\in\mathcal M$ can be brought into the subset of $\mathcal M$ having the property by an arbitrarily small perturbation that changes only the nonzero off-diagonal entries of $M$, symmetrically, and adjusts the diagonal so that row sums remain zero. Second, the perturbations of this kind that give a matrix of $\mathcal M$ without the property form a set of Lebesgue measure zero in the space $\R^{n(n-1)/2}$ of all tuples indexed by pairs $i<j$, the coordinates of pairs outside the support of $M$ being irrelevant.

Let $\mathcal M$ be the class of Laplacian matrices of connected undirected graphs with positive weights and $M=L_\Gamma(w)\in\mathcal M$. The perturbations described above produce the matrices $L_\Gamma(w+t)$, $t\in\R^E$. The first condition is \ref{ppp:dense}. For the second, $L_\Gamma(w+t)$ can belong to $\mathcal M$ only if $w+t$ has nonnegative entries. If some entry vanishes, $t$ lies in one of the finitely many hyperplanes $\{t_e=-w_e\}$; otherwise $w+t\in(0,\infty)^E$, and failure of the property means $w+t\notin\cS_\Gamma$. By \ref{ppp:null} the exceptional $t$ form a null set $N\subset\R^E$, and the corresponding set of tuples is the cylinder $N\times\R^{\binom n2-|E|}$, which is null. The same argument applies if perturbations are parametrized by the new weights $w+t$ instead of the increments $t$, since Lebesgue measure is invariant under translation. Both conditions thus hold for the conjunction of \ref{S1}--\ref{S3}. They also hold for every property implied by it, because the weights at which such a property holds contain $\cS_\Gamma$ and those at which it fails lie in the complement of $\cS_\Gamma$; this applies to ``$\lambda_2$ is simple and its eigenvector has pairwise distinct coordinates''.
\end{proof}

\begin{proof}[Proof of Corollary~\ref{cor:arithmetic}]
(i) $P_\Gamma$ is a nonzero polynomial with integer coefficients, so it cannot vanish at a point with algebraically independent coordinates. (ii) The complement of $\cS_\Gamma$ in $(0,\infty)^E$ is Lebesgue null. (iii) By the Schwartz--Zippel lemma~\cite{Schwartz,Zippel}, a nonzero polynomial of total degree $d$ over a field vanishes at a point whose coordinates are independent and uniform on a finite set $S$ with probability at most $d/|S|$. Apply this to $P_\Gamma$ over $\Q$ with $S=\{1,\ldots,N\}$. For $N=d_n+1$ the probability of failure is less than one, so a separating point of $\{1,\ldots,d_n+1\}^E$ exists.
\end{proof}

\section{Examples, symmetry and sharpness}\label{sec:examples}

\subsection{Small graphs}
The degrees~\eqref{eq:degree} are $d_2=2$, $d_3=12$, $d_4=42$, $d_5=110$ and $d_6=240$. The polynomial $P_\Gamma$ is far from minimal: it is a product of many factors, several of which never vanish on positive weights, and the examples below show much smaller exceptional sets.

\begin{example}[One edge]
For $n=2$ and weight $w$, $\kappa_1=-w$, $\kappa_2=w$ and $\delta_{12}=1$, so $P_\Gamma(w)=-w^2$. Every positive weight is separating.
\end{example}

\begin{example}[The path on three vertices]\label{ex:path3}
Let $\Gamma$ be the path $1-2-3$ with weights $a$ on $\{1,2\}$ and $b$ on $\{2,3\}$. A direct computation gives
\[
 \kappa_1=-a^2b,\quad\kappa_2=ab(a-b),\quad\kappa_3=ab^2,\qquad
 \delta_{12}=-b,\quad\delta_{13}=b-a,\quad\delta_{23}=a,
\]
and therefore
\[
 P_\Gamma(a,b)=-a^5b^5\,(a-b)^2 .
\]
The positive weights that are not separating are exactly those with $a=b$. For $a=b$ the nonconstant eigenvectors are $(1,0,-1)^T$ and $(1,-2,1)^T$: the first has a zero coordinate and the second has two equal coordinates. For $a\ne b$ the eigenvalues are $a+b\pm\sqrt{a^2-ab+b^2}$ and both nonconstant eigenvectors have nonzero, pairwise distinct coordinates.
\end{example}

\begin{example}[The triangle]\label{ex:triangle}
Let $\Gamma=K_3$ with weights $a=w_{23}$, $b=w_{13}$ and $c=w_{12}$, so that each weight is named after the opposite vertex. Then $\kappa_3=(a-b)(ab+bc+ca)$, $\delta_{12}=b-a$, and the other factors are obtained by permuting the vertices, up to sign. Altogether
\[
 P_\Gamma(a,b,c)=-(ab+bc+ca)^3\,\bigl((a-b)(b-c)(c-a)\bigr)^2 .
\]
The positive weights that are not separating are exactly the isosceles ones, with two equal weights. Indeed, a vector orthogonal to $\one$ with $x_1=x_2$ is a multiple of $(1,1,-2)^T$, and $L(1,1,-2)^T=(3b,3a,-3(a+b))^T$ is a multiple of it exactly when $a=b$; similarly $(1,-1,0)^T$ is an eigenvector exactly when $a=b$. Setting $c=0$ gives $-a^5b^5(a-b)^2$, the polynomial of the path $1-3-2$, in accordance with Lemma~\ref{lem:integrality}(ii) and Example~\ref{ex:path3}.
\end{example}

\subsection{Symmetric weightings are exceptional}
An \emph{automorphism} of the weighted graph $(\Gamma,w)$ is a permutation $\sigma$ of $V$ that maps edges to edges and preserves their weights; equivalently, its permutation matrix commutes with $L_\Gamma(w)$.

\begin{proposition}\label{prop:symmetry}
Let $\Gamma$ be a connected graph on $n\ge3$ vertices and $w\in\cS_\Gamma$. Then the weighted graph $(\Gamma,w)$ has no automorphism other than the identity. In particular, if $\Gamma$ has a nontrivial automorphism, the constant weights are not separating.
\end{proposition}
\begin{proof}
Let $\Pi$ be the permutation matrix of an automorphism $\sigma$, so that $\Pi L=L\Pi$ for $L=L_\Gamma(w)$. If $x$ is a unit eigenvector, then $\Pi x$ is a unit eigenvector for the same simple eigenvalue, so $\Pi x=\pm x$. Suppose $\Pi x=x$ for some eigenvector $x\perp\one$. Then $x$ is constant on the cycles of $\sigma$, and \ref{S3} forces $\sigma$ to be the identity. Otherwise $\Pi x=-x$ for every eigenvector orthogonal to $\one$, and $\Pi\one=\one$, so $\Pi=\tfrac2n\one\one^T-I$. The off-diagonal entries of this matrix equal $2/n$, which is neither $0$ nor $1$ when $n\ge3$, so it is not a permutation matrix.
\end{proof}

For $n=2$ the transposition is an automorphism and every weight is separating, so the hypothesis $n\ge3$ is needed. Every graph on at most five vertices has a nontrivial automorphism; accordingly $P_\Gamma$ vanishes at the constant weights for all connected graphs with $3\le n\le5$. This is the phenomenon behind the repeated or vanishing Fiedler entries of small unweighted graphs.

\subsection{Exceptional weightings without symmetry}
In Examples~\ref{ex:path3} and~\ref{ex:triangle} the exceptional weightings are exactly the symmetric ones. This is not the general situation. The next two examples exhibit Fiedler vectors with a zero coordinate, and with two equal coordinates, on weighted graphs whose only automorphism is the identity.

\begin{example}[A path with a vanishing Fiedler coordinate]\label{ex:path4}
Let $\Gamma$ be the path $1-2-3-4$ with weights $(w_{12},w_{23},w_{34})=(1,3,2)$. Then
\[
 L=\begin{pmatrix}1&-1&0&0\\-1&4&-3&0\\0&-3&5&-2\\0&0&-2&2\end{pmatrix},\qquad
 L\begin{pmatrix}-3\\0\\1\\2\end{pmatrix}=\begin{pmatrix}-3\\0\\1\\2\end{pmatrix},
\]
and the characteristic polynomial is $z(z-1)(z-3)(z-8)$. Thus $\lambda_2=1$ is simple and the Fiedler vector $(-3,0,1,2)^T$ vanishes at vertex~$2$. The weights are not palindromic, so the reflection of the path is not an automorphism.
\end{example}

\begin{example}[A cycle with two equal Fiedler coordinates]\label{ex:cycle4}
Let $\Gamma$ be the cycle $1-2-3-4-1$ with weights $(w_{12},w_{23},w_{34},w_{41})=(8,15,35,12)$. Then
\[
 L=\begin{pmatrix}20&-8&0&-12\\-8&23&-15&0\\0&-15&50&-35\\-12&0&-35&47\end{pmatrix},\qquad
 L\begin{pmatrix}-4\\1\\2\\1\end{pmatrix}=25\begin{pmatrix}-4\\1\\2\\1\end{pmatrix},
\]
and the characteristic polynomial is $z(z-25)(z^2-115z+2448)$. The quadratic factor is positive at $z=25$ and has its minimum at $z=115/2>25$, so both of its roots exceed $25$. Hence $\lambda_2=25$ is simple, and the Fiedler vector $(-4,1,2,1)^T$ takes the same value at the vertices $2$ and $4$. The four weights are distinct, and no nonidentity automorphism of the $4$-cycle fixes all four edges, so the weighted graph has trivial automorphism group.
\end{example}

\subsection{Paths}
On a path the conclusion of the conjecture holds for every positive weighting.

\begin{proposition}[Paths]\label{prop:paths}
Let $\Gamma$ be the path $1-2-\cdots-n$, $n\ge2$, with positive weights $w_i=w_{i,i+1}$. Every eigenvalue of $L_\Gamma(w)$ is simple, and every Fiedler vector is strictly monotone along the path; in particular its coordinates are pairwise distinct.
\end{proposition}
\begin{proof}
Write $L=L_\Gamma(w)$. If $Lx=\lambda x$, the first row gives $w_1(x_1-x_2)=\lambda x_1$, and row $k$, for $1<k<n$, gives $w_k(x_k-x_{k+1})=\lambda x_k-w_{k-1}(x_k-x_{k-1})$. Since every $w_k$ is positive, $x_1$ determines $x_2,\ldots,x_n$ recursively. Hence every eigenspace has dimension one, so every eigenvalue is simple, and every eigenvector has $x_1\ne0$.

Let $x$ be a Fiedler vector and put $d_i=x_{i+1}-x_i$ and $S_i=x_1+\cdots+x_i$ for $1\le i<n$. For every $z\perp\one$,
\[
 \|z\|^2=\frac1n\sum_{1\le i<j\le n}(z_j-z_i)^2 ,
\]
since the right side equals $\|z\|^2-\frac1n(\one^Tz)^2$. Let $y\perp\one$ be the vector with $y_{i+1}-y_i=|d_i|$. By Lemma~\ref{lem:psd}, $y^TLy=\sum_iw_id_i^2=x^TLx=\lambda_2\|x\|^2$, while
\[
 \|y\|^2=\frac1n\sum_{i<j}\Bigl(\sum_{k=i}^{j-1}|d_k|\Bigr)^2
 \ \ge\ \frac1n\sum_{i<j}\Bigl(\sum_{k=i}^{j-1}d_k\Bigr)^2=\|x\|^2 .
\]
On $\one^\perp$ the Rayleigh quotient of $L$ is at least $\lambda_2$ by the Courant--Fischer theorem~\cite[Section~4.2]{HornJohnson}, and $y\ne0$ because $x$ is not constant. Hence $\lambda_2\|y\|^2\le y^TLy=\lambda_2\|x\|^2$, and equality holds throughout. The term $(i,j)=(1,n)$ then gives $|\sum_kd_k|=\sum_k|d_k|$, so the nonzero $d_k$ share one sign; replacing $x$ by $-x$, we may assume that $x$ is nondecreasing.

Suppose that $d_i=0$ for some $i<n$. Lemma~\ref{lem:cut}, applied to $D=\{i+1,\ldots,n\}$ and the edge $\{i,i+1\}$, gives $w_id_i=\lambda_2(x_{i+1}+\cdots+x_n)=-\lambda_2S_i$, because $x\perp\one$; hence $S_i=0$. Monotonicity gives $0=S_i\le ix_i$ and $0=x_{i+1}+\cdots+x_n\ge(n-i)x_{i+1}=(n-i)x_i$, so $x_i=0$. Then $x_1\le\cdots\le x_i=0$ and $S_i=0$ force $x_1=0$, which contradicts the first paragraph. Therefore every $d_i$ is nonzero, and $x$ is strictly increasing.
\end{proof}

The Fiedler vector of a path can nevertheless vanish at a vertex, and the other eigenvectors can repeat a coordinate, as Examples~\ref{ex:path3} and~\ref{ex:path4} show; so \ref{S2} and \ref{S3} may fail on paths.

\subsection{Scope of the hypotheses}

\begin{remark}[Connectedness]
If $\Gamma$ is disconnected, $P_\Gamma=0$ by Lemma~\ref{lem:integrality}(iii), and no weight vector, positive or not, satisfies~\ref{S1}. Connectedness is thus necessary and sufficient for $P_\Gamma\ne0$.
\end{remark}

\begin{remark}[Real weights]
Proposition~\ref{prop:certificate} does not use positivity. Hence, for a connected graph, the matrix $L_\Gamma(w)$ satisfies \ref{S1}--\ref{S3} for every real weight vector $w$ outside the proper algebraic set $Z(P_\Gamma)\subset\R^E$, for instance for generic signed weights. Positivity is used only to identify $\lambda_2$ as the smallest nonzero eigenvalue.
\end{remark}

\section{Remarks and questions}\label{sec:questions}

\begin{remark}[Why a first-order expansion is unavoidable]
In the induction of Section~\ref{sec:trees} the matrix $L_0=L\oplus0$ is not separating, and $P_{T'}(w,0)=0$. The nonvanishing of $P_{T'}(w,\varepsilon)$ for small $\varepsilon>0$ is therefore a statement about the lowest-order term of this polynomial in $\varepsilon$. Lemma~\ref{lem:neweigenvector} computes that term for the factors $\delta_{ij}$ with $i,j\in V$: by~\eqref{eq:deltaformula}, each of them vanishes to order exactly one in $\mu$ when $h_i\ne h_j$.
\end{remark}

\begin{remark}[Spectral orderings]
For separating weights, sorting the vertices by the Fiedler vector gives a strict total order, unique up to reversal, and the sign of the Fiedler vector splits $V$ into two nonempty sets without ties. By Fiedler's theorem~\cite{Fiedler75} each of the two sets induces a connected subgraph. Theorem~\ref{thm:main} says that, for a fixed graph, ties in these constructions occur only on a null algebraic set of weights.
\end{remark}

\begin{question}[The exceptional hypersurface]\label{q:degree}
Let $\Gamma$ be connected. What is the smallest degree of a polynomial whose zero set in $(0,\infty)^E$ is the complement of $\cS_\Gamma$? For the path on three vertices the answer is $1$, while $d_3=12$. Which irreducible components does the exceptional set have, and which of them correspond to a zero coordinate, to two equal coordinates, or to a multiple eigenvalue?
\end{question}

\begin{question}[Quantitative separation]
For random weights with a bounded density, how small can $\min_{i\ne j}|\phi(i)-\phi(j)|$ be for the unit Fiedler vector $\phi$? A polynomial lower bound in $n$ with high probability would make the spectral order numerically stable.
\end{question}

\begin{question}[Other operators]
Do the analogues of Theorem~\ref{thm:main} hold for normalized Laplacians, for adjacency matrices with a fixed support, and for directed graphs with a spanning diverging tree, the setting in which simple spectrum is generic by~\cite[Theorem~6.6]{PPP}? The Krylov certificate of Section~\ref{sec:certificate} adapts to any family of symmetric matrices; the construction of a witness is the part that depends on the operator.
\end{question}


\begin{thebibliography}{99}
\bibitem{Fiedler73} M. Fiedler. Algebraic connectivity of graphs. \emph{Czechoslovak Mathematical Journal} 23 (1973), 298--305.
\bibitem{Fiedler75} M. Fiedler. A property of eigenvectors of nonnegative symmetric matrices and its application to graph theory. \emph{Czechoslovak Mathematical Journal} 25 (1975), 619--633.
\bibitem{ABH} J. E. Atkins, E. G. Boman, and B. Hendrickson. A spectral algorithm for seriation and the consecutive ones problem. \emph{SIAM Journal on Computing} 28 (1998), 297--310.
\bibitem{PPP} C. Poignard, T. Pereira, and J. P. Pade. Spectra of Laplacian matrices of weighted graphs: structural genericity properties. \emph{SIAM Journal on Applied Mathematics} 78(1) (2018), 372--394. \url{https://doi.org/10.1137/17M1124474}. arXiv:1704.01677.
\bibitem{Kalman} R. E. Kalman. On the general theory of control systems. \emph{Proceedings of the First IFAC Congress} (1960).
\bibitem{Uhlenbeck} K. Uhlenbeck. Generic properties of eigenfunctions. \emph{American Journal of Mathematics} 98 (1976), 1059--1078.
\bibitem{Godsil} C. Godsil. Controllable subsets in graphs. \emph{Annals of Combinatorics} 16 (2012), 733--744.
\bibitem{ORT} S. O'Rourke and B. Touri. On a conjecture of Godsil concerning controllable random graphs. \emph{SIAM Journal on Control and Optimization} 54 (2016), 3347--3378.
\bibitem{TaoVu} T. Tao and V. Vu. Random matrices have simple spectrum. \emph{Combinatorica} 37 (2017), 539--553.
\bibitem{BLS} T. B{\i}y{\i}ko\u{g}lu, J. Leydold, and P. F. Stadler. \emph{Laplacian Eigenvectors of Graphs: Perron--Frobenius and Faber--Krahn Type Theorems}. Lecture Notes in Mathematics 1915, Springer, 2007.
\bibitem{HornJohnson} R. A. Horn and C. R. Johnson. \emph{Matrix Analysis}. Second edition, Cambridge University Press, 2013.
\bibitem{Schwartz} J. T. Schwartz. Fast probabilistic algorithms for verification of polynomial identities. \emph{Journal of the ACM} 27 (1980), 701--717.
\bibitem{Zippel} R. Zippel. Probabilistic algorithms for sparse polynomials. In \emph{Symbolic and Algebraic Computation (EUROSAM '79)}, Lecture Notes in Computer Science 72, Springer, 1979, pp.~216--226.
\end{thebibliography}
\end{document}